\documentstyle[12pt,psfig]{article}
\newcommand{\copyleft}{
GNU FDL\thanks{
Copyright (C) 1994 Peter G. Doyle.
Permission is granted to copy, distribute and/or modify this document
under the terms of the GNU Free Documentation License, 
as published by the Free Software Foundation;
with no Invariant Sections, no Front-Cover Texts, and no Back-Cover Texts.
}}

\title{Division by three}
\author{
Peter G. Doyle
\and
John Horton Conway
\thanks{John Conway collaborated on the research reported here, and has
been listed as an author of this work since it was first distributed in 1994.
But he has never approved of this exposition, which he regards
as full of `fluff'.}
}
\date{Version dated 1994
\\ \copyleft}

\newcommand{\fig}[3]{
\begin{figure}
\psfig{figure=figures/#1.ps,width=370pt}
\caption{#3}
\label{#2}
\end{figure}
}

\newcommand{\qed}{\clubsuit}
\newcommand{\dispqed}{\;\;\qed}
\newcommand{\union}{\cup}
\newcommand{\cross}{\times}
\newcommand{\plus}{+}
\newcommand{\implies}{\Longrightarrow}

\newcommand{\seteq}{\asymp}
\newcommand{\setleq}{\preceq}

\newcommand{\setgeq}{\succeq}
\newcommand{\setll}{\ll}
\newcommand{\setgg}{\gg}
\newcommand{\degree}{^\circ}

\begin{document}

\maketitle

\begin{abstract}
We prove without appeal to the Axiom of Choice that for any sets $A$ and $B$,
if there is a one-to-one correspondence
between $3 \cross A$ and $3 \cross B$
then there is a one-to-one correspondence between $A$ and $B$.
The first such proof, due to Lindenbaum,
was announced by Lindenbaum and Tarski in 1926,
and subsequently `lost'; Tarski published an alternative proof
in 1949.
We argue that the proof presented here follows Lindenbaum's
original.

\vspace{0.5cm}
{\bf AMS Classification numbers}\quad
03E10 (Primary); 03E25 (Secondary).
\end{abstract}

\section{Introduction}
In this paper we show that it is possible to divide by three.
Specifically,
we prove that for any sets $A$ and $B$,
if there is a one-to-one correspondence
between $3 \cross A$ and $3 \cross B$
then there is a one-to-one correspondence between $A$ and $B$.
Here $3$ denotes the standard three-element set $\{0,1,2\}$
and $\cross$ denotes Cartesian product.
Writing $A \seteq B$ to indicate that there is a one-to-one correspondence
between $A$ and $B$,
our assertion is that
\[
3 \cross A \seteq 3 \cross B \implies A \seteq B
.
\]

This assertion is easy to prove if we allow ourselves to use
the axiom of choice,
which states (roughly) that it is always possible to choose one element
simultaneously from each one of an infinite collection of non-empty sets.
If we want to avoid using the axiom of choice---and we do---then
the problem of dividing by three becomes much more difficult,
and while the problem was discussed by Bernstein in 1901,
and a solution announced by Lindenbaum and Tarski in 1926,
the earliest published solution is Tarski's solution of 1949.
We'll say more about the history of the problem later on.

Now when we say that we want to avoid using the axiom of choice,
it is not merely, or even mainly, because we think the axiom is probably
not `true'.
The real reason for avoiding the axiom of choice is that it is not
{\em definite}.
(See Jech
\cite{jech:setTheory}.)
To prove that $3 \cross A \seteq 3 \cross B$ means somehow or other
to show that from a one-to-one correspondence between $3 \cross A$
and $3 \cross B$,
you can produce a one-to-one correspondence between $A$ and $B$.
If the proof requires the axiom of choice,
then the procedure will involve at some point or other a
`subroutine call' to the axiom of choice,
where we pass in a collection of non-empty sets,
and receive as output a `choice function' selecting one element from
each of these sets,
and the correspondence between $A$ and $B$ that we end up with
will depend on which particular choice function was `returned'.
Thus the correspondence will not be uniquely defined.
Conversely,
if we avoid the axiom of choice,
then the procedure will produce a well-defined correspondence between
$A$ and $B$.

By the same token,
if we make sure that all our constructions are well-defined,
then we can make sure that we haven't used the axiom of choice.
This means that we will be able to do set theory in the usual
free-and-easy way,
without sticking close to the axioms of ZF or any other formal
system for doing set theory.
Note that in order to steer clear of the axiom of choice it is not
necessary to avoid making any choices at all (what would that mean?),
or even to avoid making any apparently arbitrary choices
(e.g. choices that break some natural symmetry of the problem
we're trying to solve), but simply to make sure
that when there is a choice to be made,
we describe what to do in such a way that there is one and only one
way to do it.

To sum up, our attitude towards the axiom of choice is this:
To avoid the axiom of choice is the same thing as to make sure
that all constructions are uniquely defined.
It is the problem of finding a uniquely defined procedure for dividing by
three that we are really interested in.
This problem has connections to constructive (finite) combinatorics
and to parallel processing,
but never mind all that:
The problem is
easy to state,
hard to solve, and just plain fun.
Thus we hope that that this paper will appeal to you even if you
can't muster much anxiety about the axiom of choice.

You're probably wondering about dividing by two.
The reason this paper isn't called `Division by two'
is that dividing by two is significantly easier than dividing
by three,
though still not trivial.
Three seems to be the crucial case:
the methods we will discuss for dividing by three extend
to allow you to divide by 5
(you can already divide by 4 if you can divide by 2:
just do it twice),
by 7,
and in general by any finite $n$.
What's more, you can use the same method to show that if
$3 \cross A \seteq 2 \cross B$
then there is a set $C$ such that
$A \seteq 2 \cross C$ and $B \seteq 3 \cross C$,
and so on.
(See Bachmann
\cite{bachmann:transfinite}.)
But here we will only explicitly discuss
the cases of division by two and by three.

\section{What makes this problem difficult?}

In order to show that $3 \cross A \seteq 3 \cross B \implies A \seteq B$,
what we need to show is that given a one-to-one correspondence between
$3 \cross A$ and $3 \cross B$,
we can come up with a one-to-one correspondence between $A$ and $B$.
What makes this difficult is that this input correspondence is all we
have to work with:
We're not allowed to ask someone to manufacture a well-ordering of
$A$
(or anything else, for that matter)
out of thin air.

Now it is remarkable that it should be hard to divide by three,
so let's think about why it seems at first glance to be a trivial
problem.

The assertion
$3 \cross A \seteq 3 \cross B \implies A \seteq B$
about sets has an aritmetical counterpart,
namely that if $a$ and $b$ are natural numbers then
$3a=3b \implies a=b$.
This simple arithmetical statement has a simple proof
from the axioms of arithmetic
(never mind {\em which} axioms of arithmetic---any reasonable
set of axioms will do),
and it is would be surprising indeed if this proof
could not be converted into a proof that
$3 \cross A \seteq 3 \cross B \implies A \seteq B$
in the case where $A$ (and hence also $B$) is finite.

To make this conversion,
we need to consider briefly the notion of a 
{\em cardinal number} ({\em cardinal} for short).
A cardinal number is defined
to be an equivalence class of sets under the relation $\seteq$,
and if we make appropriate definitions for sum and product,
then our statement
$3 \cross A \seteq 3 \cross B \implies A \seteq B$
translates into the statement
$3a=3b \implies a=b$.
The {\em finite cardinals}
$[0],[1],[2],[3],\ldots$
are the $\seteq$-equivalence classes of the standard finite sets
$0=\emptyset,1=\{0\},2=\{0,1\},3=\{0,1,2\},\ldots$.
It is not difficult to prove that the finite cardinals satisfy the 
axioms for arithmetic,
and combining a proof that the finite cardinals satisfy the axioms of arithmetic
with a proof from the axioms of arithmetic that $3a=3b \implies a=b$
yeilds a proof that 
$3 \cross A \seteq 3 \cross B \implies A \seteq B$
as long as $A$ (and hence also $B$) is a finite set.

To deal with the infinite case,
we might be tempted to argue
that if $A$ is infinite then $3 \cross A \seteq A$,
and so when $3 \cross A \seteq 3 \cross B$ and $A$ (and hence also $B$)
is infinite,
$A \seteq 3 \cross A \seteq 3 \cross B \seteq B$.
The problem with this argument is that the assertion that
$3 \cross A \seteq A$ for infinite $A$ depends on the axiom of choice.
So much for that idea.

Another idea would be to adapt to the general case
the solution that we got (in principle)
in the finite case.
If we try to do this,
we will find that the proof in the finite case is not really as good
as we might have imagined,
for while it will give a procedure for converting a one-to-one
correspondence between $3 \cross A$ and $3 \cross B$
into a one-to-one correspondence between $A$ and $B$,
when you try to use this procedure you will be called upon to
input not only the requisite
correspondence between $3 \cross A$ and $3 \cross B$,
but also a bijection between $A$ and $\{0,1,2,\ldots,n-1\}$,
and the correspondence between $A$ and $B$ that you receive as output
will depend on this second input as well as the first.
Because it depends on an extraneous input,
namely a particular ordering of the set $A$,
the procedure will not be `canonical', or `equivariant',
and it will not extend to the infinite case because then
we will not have any ordering of the set $A$ to call upon.
As you may know,
the axiom of choice is equivalent to the statement that
`any set can be well-ordered',
and a well-ordering of $A$ is exactly what we would need to
push this method of dividing by three through in the infinite case.
So this isn't going to work either.

But this is getting complicated, whereas the problem is really
very simple:
Turn the correspondence between $3 \cross A$ and $3 \cross B$
into a correspondence between $A$ and $B$.

\section{History}
A proof that it is possible to divide by two was presented by
Bernstein in his Inaugural Dissertation of 1901,
which appeared in Mathematische Annallen in 1905
\cite{bernstein:mengenlehre};
Bernstein also indicated how to extend his results to division by any
finite $n$,
but we are not aware of anyone other than Bernstein himself
who ever claimed to understand this argument.
In 1922 Sierpi\'{n}ski
\cite{sierpinski:two}
published a simpler proof of division by two,
and he worked hard to extend his method to division by three, but never
succeeded.

In 1926 Lindenbaum and Tarski announced,
in an infamous paper
\cite{lindenbaumTarski:ensembles}
that contained statements (without proof) of
144 theorems of set theory,
that Lindenbaum had found a proof
of division by three.
Their failure to give any hint of a proof must have frustrated Sierpi\'{n}ski,
for it appears that twenty years later he still did not know how to
divide by three
(cf.
\cite{sierpinski:leq}).
Finally, in 1949, in a paper
`dedicated to Professor Waclaw Sierpi\'{n}ski in celebration of his forty years
as teacher and scholar',
Tarski
\cite{tarski:cancellation}
published a proof.
In this paper,
Tarski explained that unfortunately he couldn't remember how Lindenbaum's
proof had gone,
except that it involved an argument like the one
Sierpi\'{n}ski had used in dividing by two,
and another lemma, due to Tarski, which we will describe below.
Instead of Lindenbaum's proof,
he gave another.

Now when we began the investigations reported on here,
we were aware that there was a proof in Tarski's paper,
and Conway had even pored over it at one time or another
without achieving enlightenment.
The problem was closely related to the kind of question John had looked
at in his thesis,
and it was also related to work
that Doyle had done in the field of bijective
combinatorics.
So we decided that we were going to figure out what the heck was going on.
Without too much trouble we figured out how to divide by two.
Our solution turned out to be substantially equivalent to that of
Sierpi\'{n}ski,
though the terms in which we will describe it below will not much resemble
Sierpi\'{n}ski's.
We tried and tried and tried to adapt the method to the case of dividing by
three, but we kept getting stuck at the same point in the argument.
So finally we decided to look at Tarski's paper,
and we saw that the lemma Tarski said Lindenbaum had used was precisely
what we needed to get past the point we were stuck on!
So now we had a proof of division by three that combined an argument like
that Sierpi\'{n}ski used in dividing by two with an appeal to Tarski's lemma,
and we figured we must have hit upon an argument very much like that of
Lindenbaum's.
This is the solution we will describe here:
Lindenbaum's argument, after 62 years.

\section{The Cantor-Schr\"{o}der-Bernstein theorem}

Probably the best way to see what it means to do set theory without the
axiom of choice
is to look at the proof of the Cantor-Schr\"{o}der-Bernstein theorem,
which is the prototypical example.
Let us write $A \setleq B$ when there is an injection from $A$ into $B$.

{\bf Cantor-Schr\"{o}der-Bernstein Theorem.}
If $A \setleq B$ and $B \setleq A$ then $A \seteq B$.

{\bf Proof.} We want to show that
given injections $f:A \to B$ and $g:B \to A$
we can determine a one-to-one correspondence between $A$ and $B$.
We can and will assume that $A$ and $B$ are disjoint.
Here's how it goes.
We visualize the set $A$ as a collection of blue dots,
and the set $B$ as a collection of red dots.
(See Figure \ref{fig:csbgraph}.)
\fig{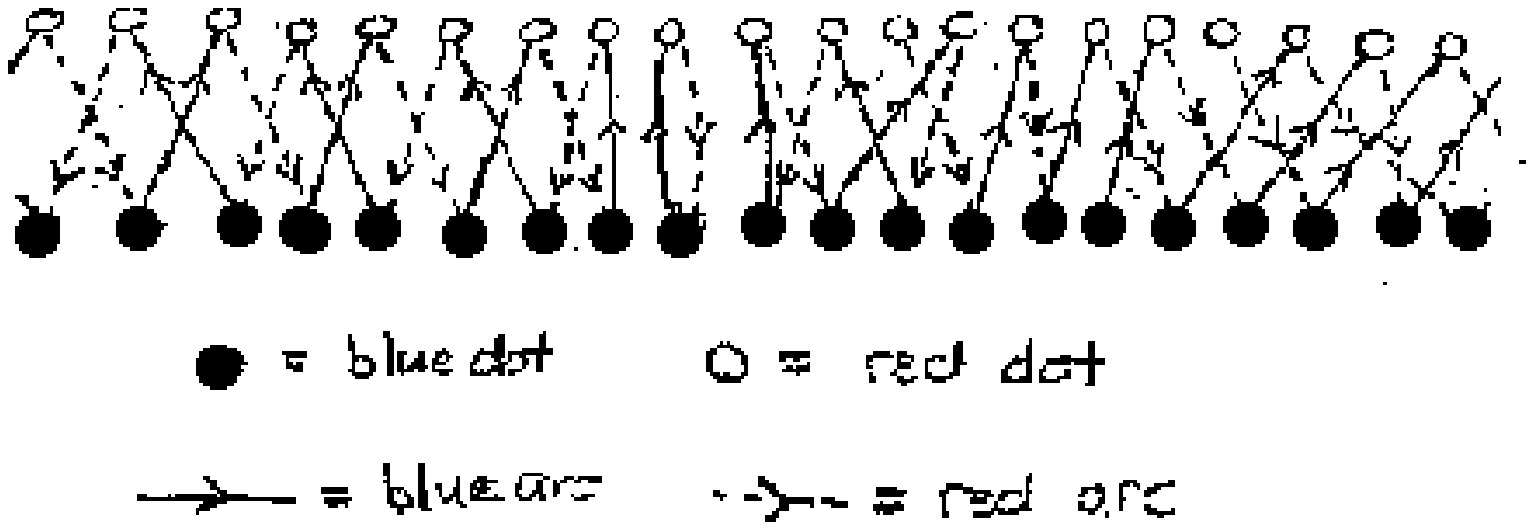}{fig:csbgraph}
{
The blue dots represent elements of $A$;
the red dots represent elements of $B$.
The blue arcs connect pairs $(x,f(x))$;
the red arcs connect pairs $(y,g(y))$.
}
We visualize the injection $f$ as a collection of blue directed arcs
connecting each element $x \in A$ to its image $f(x) \in B$.
Similarly, we visualize $g$ as a collection of red directed arcs.
If we put in both the blue and the red arcs,
we get a directed graph where
every vertex has one arc going out and at most one arc coming in.

Such a graph decomposes into a union of connected components,
each of which is either a finite directed cycle,
a doubly-infinite path,
or a singly-infinite path.
As you go along one of these paths or cycles, the vertices you encounter
belong alternately to $A$ and $B$.
In the case of a cycle or a doubly-infinite path,
the blue arcs define a one-to-one correspondence between the
blue vertices of the component and the red vertices.
In the case of a singly-infinite path,
the blue edges will still determine a one-to-one correspondence between
the blue and red vertices of the path if the path begins with a blue vertex,
but not if the path begins with a red vertex.
However in this latter case
(the `Sadie Hawkins case')
we can take the red edges instead.
Thus we can pair up the vertices of $A$ and $B$ along each connected
component,
and the union of these correspondences determines a one-to-one
correspondence between $A$ and $B$.
(See Figure \ref{fig:csbmatch}.)
$\qed$
\fig{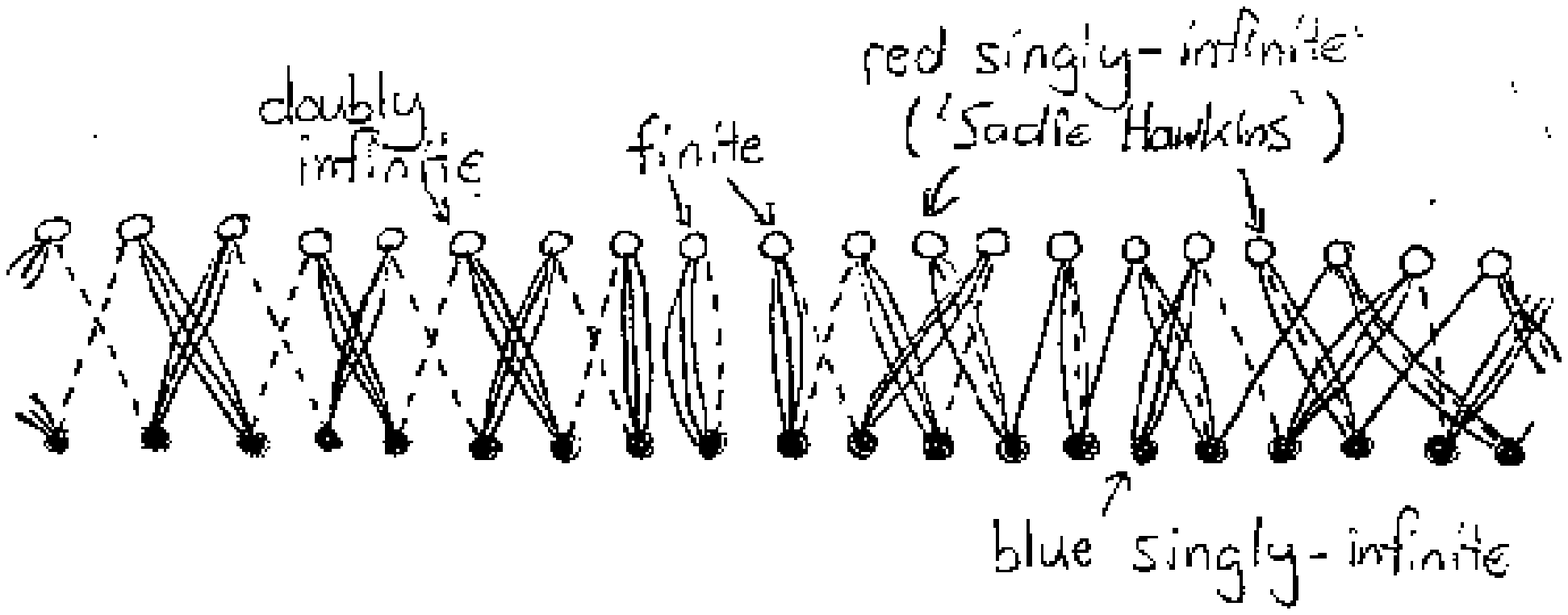}{fig:csbmatch}
{
The matching follows the blue edges except on Sadie Hawkins components,
where it follows the red edges.
}

The foregoing argument is a good illustration of how to do set theory
without the axiom of choice.
(Or rather, it is a good illustration of one way to do set theory
without the axiom of choice---there are other ways as well.)
Beginning with $f$ and $g$,
we describe a graph built using the information contained in $f$ and $g$,
and we determine our one-to-one correspondence
by fooling around on the graph.

Note that while the fooling around on the graph is definite,
it is not `effectively computable' or anything like that.
For instance, suppose you want to know what element $y$ of $B$ corresponds
to a particular element $x$ of $A$.
Most likely $y=f(x)$,
because this is the element of $B$ that is joined to $x$ by a blue arc
(oriented out of $x$).
However, if $x$ happens to belong to a singly-infinite path
that begins with an element of $B$,
then the corresponding element of $B$ is the unique $y$ such that
$g^{-1}(\{x\})=\{y\}$,
because this is the element of $B$ that is joined to $x$ by a red arc
(oriented into $x$, in this case).
Now suppose that $f$ and $g$ and their inverses are functions that
we are able to compute,
so that for example for any $x$ in $A$ we can determine $f(x)$,
and we can determine
if there is any $y$ such that $x = g(y)$,
and if so we can determine what $y$ is.
To compute the element corresponding to a given $x \in A$,
we need to know whether
the path that starts at $x$ and goes backwards along the edges of the
directed graph terminates at an element of $B$.
When $A$ and $B$ are infinite,
this is not something we can expect to determine in finite time,
because while we can tell what the answer is if the path 
terminates (or connects back up to the starting point $x$),
if it hasn't terminated after a long time that doesn't mean that it will
never terminate.
In order to consider our matching as the result of some kind of
algorithm, we would need to allow a certain amount of infiniteness into
our computations:
We would need to be able to say,
for example,
``Follow that path as far as it goes,
even if it continues on forever,
and tell me whether it ever terminates or not.''
If we are allowed this kind of infinite computation,
then the correspondence here becomes computable,
and in fact the whole theory of what can be explicitly determined
comes to coincide with the theory of what you can compute.

But all this bother about the extent to which the correspondence we've
defined is computable is not relevant to our present purpose,
which is to be sure to avoid an appeal to the axiom of choice.
The description we've given determines uniquely a
one-to-one correspondence between $A$ and $B$,
and that's all we need to make sure of.

The Cantor-Schr\"{o}der-Bernstein theorem shows that
$\setleq$ induces a partial ordering $\leq$ of cardinals if we
say that $a \leq b$ whenever $a=[A]$, $b=[B]$, and $A \setleq B$.
However, without the axiom of choice,
this relation is not necessarily a total ordering,
because if you live in a world of sets where the axiom of choice does not hold,
it is not necessarily true that either $A \setleq B$ or $B \setleq A$.
The failure of $\setleq$ to induce a total order on cardinals
is a clue that without the axiom of choice,
the whole notion of cardinal number is not all that
natural or useful.
Another indication is this:
There is another definition (due to Dedekind)
of what it means for a set (and its cardinal)
to be infinite;
we will discuss this alternative notion of infinite sets later on,
when we discuss subtraction.
In the absense of the axiom of choice this second definition
does not necessarily agree with the definition above
that a set is
finite if there is a bijection between it and one of the standard finite
sets $\{0,1,2,\ldots,n-1\}$.
For these reasons, we will avoid as much as possible any mention of
cardinals in what follows.

\section{Division by two}

We are given a bijection
$f:\{0,1\} \cross A \rightarrow \{0,1\} \cross B$,
and we want to produce a bijection
$g:A \rightarrow B$.

The key to our construction will be to have a clear mental picture
of the original bijection $f$.
To this end, we imagine that the set $A$ consists of
a whole bunch of blue arrows---maybe they're made of wood---and
the set $B$ consists of a whole bunch of red arrows.
Each arrow has a distinguishable head and tail.
If we label the tail of each arrow 0 and the head 1,
then we can identify the set of `ends of blue arrows' with
the set $\{0,1\} \cross A$, and similarly for the red arrows.
The bijection $f$ now provides a one-to-one correspondence between
ends of blue arrows and ends of red arrows, so that to
each end of each arrow there corresponds a uniquely determined end of
an arrow of the opposite color.
Somehow,
we must convert this matching between ends of arrows into a
matching between the arrows themselves.

Now imagine that the matching between ends of arrows
constitutes  the {\em Instructions} for assembling
a toy that you are supposed to  build using the arrows supplied.
Each end of a blue arrow is to be connected to the corresponding
end of a red arrow.
When they're all hooked together,
some arrows will meet head to head, some tail to tail, some head to tail.
Each arrow will belong to a unique string of arrows,
either a finite string joined around in a sort of necklace,
or an infinite string of arrows stretching off in both directions.
On each string, the arrows alternate in color.
(Don't worry about whether the strings are tangled together or knotted;
this kind of thing isn't of any concern to us.)

Here's an example.
Suppose that $A = \{a,b,c\}$ and
$B = \{x,y,z\}$, and that the correspondence $f$ is as illustrated
in Figure \ref{fig:matching}.
\fig{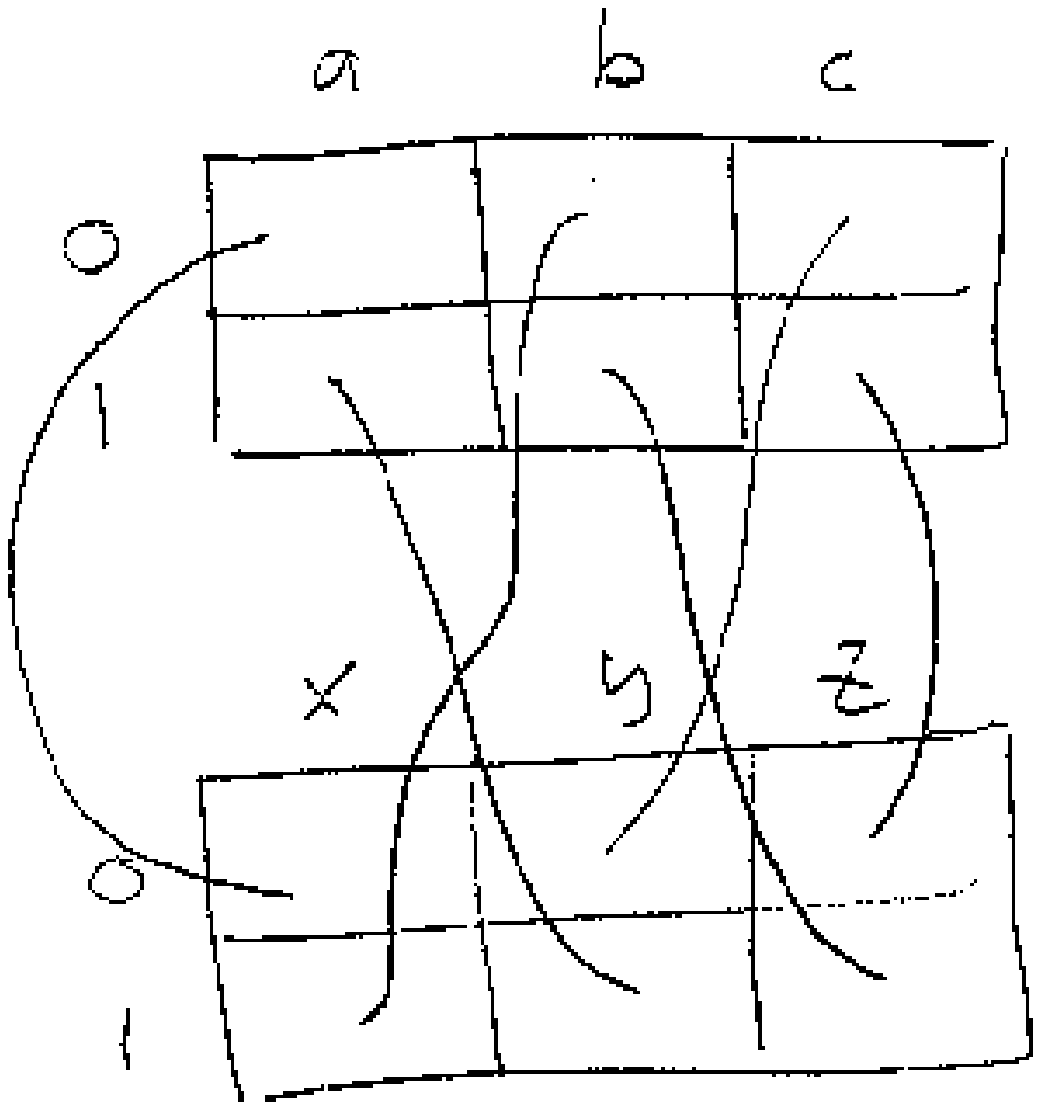}{fig:matching}
{
A matching between $\{0,1\} \cross \{a,b,c\}$
and $\{0,1\} \cross \{x,y,z\}$.
}
The head of the blue arrow $a$ is represented by the pair
$(1,a)$, and according to the diagram this is matched to the pair
$(1,y)$, which represents the head of $y$,
so we connect together the heads of $a$ and $y$.
Now since $(0,y)$ is matched to $(0,c)$,
we hook the tail of $y$ to the tail of $c$,
and so forth.
When all of the connections are made,
the six arrows join together to form a necklace,
as illustrated in Figure \ref{fig:necklace}.
\fig{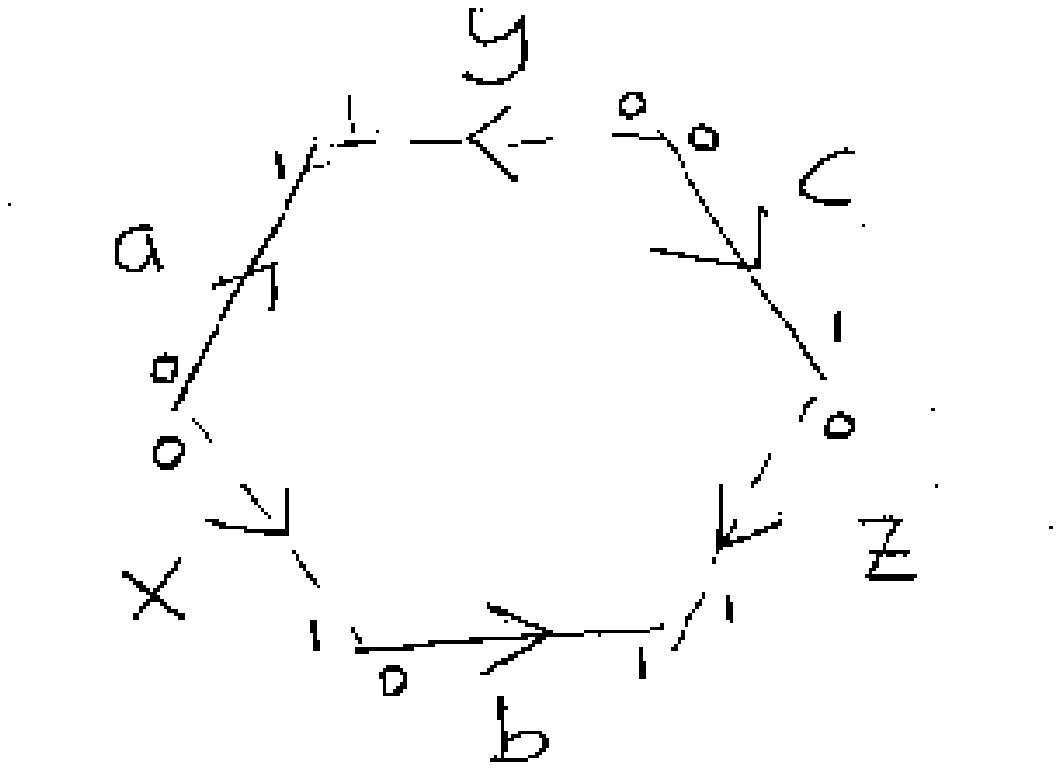}{fig:necklace}
{
The arrows are joined as indicated by the correspondence
in the previous figure.
Note that whereas in earlier figures arcs merely connected together
the objects of primary interest,
here it is the arrows themselves
(they're {\em arrows} now, not arcs) that are primary.
}

Now since the arrows alternate in color,
it is intuitively clear that on any string there are just as many red
parentheses as blue parentheses,
and we ought to be able to match up the arrows separately on each string,
so that each arrow gets matched an arrow of the opposite color
on the same string.
All we have to do is tell how to do it.

Imagine for a moment
that the strings of arrows represent streets---circular
drives, in the case of necklaces, and long boulevards in the case of
infinite strings.
If the arrows are good, straight, {\em American} arrows,
it is very natural for each arrow to dream of marrying the arrow next door.
The only difficulty is that there are two arrows next door,
because there are two next doors.
Of course since each arrow has a direction associated to it,
it might be that they
have in mind the arrow that their head is connected to, but this is
going to cause all kinds of conflicts.
For instance, in the example above, both $b$ and $c$ are going to want
to marry $z$.

If only all the streets were one-way streets,
then the blue arrows could resolve to marry the red arrow just up the street,
and the red arrows could resolve to marry the blue arrow just down the street,
and everything would be dandy.
If only someone would go around and decide
somehow or other for each street
which way the traffic should go!
But this is exactly the kind of thing we're not allowed to ask for,
if we want to steer clear of the axiom of choice.
Unless we can describe a rule for determining the direction along a street,
we're out of luck.

Now if you look back at Figure \ref{fig:necklace},
you'll notice that the arrows that are shown there
are not real arrows with arrowheads and feathers, but mathematicians'
arrows with the direction indicated by a sort of angle-bracket gadget
in the middle.  If we round the angle-brackets into parentheses, we
get the picture shown in Figure \ref{fig:parens}.
\fig{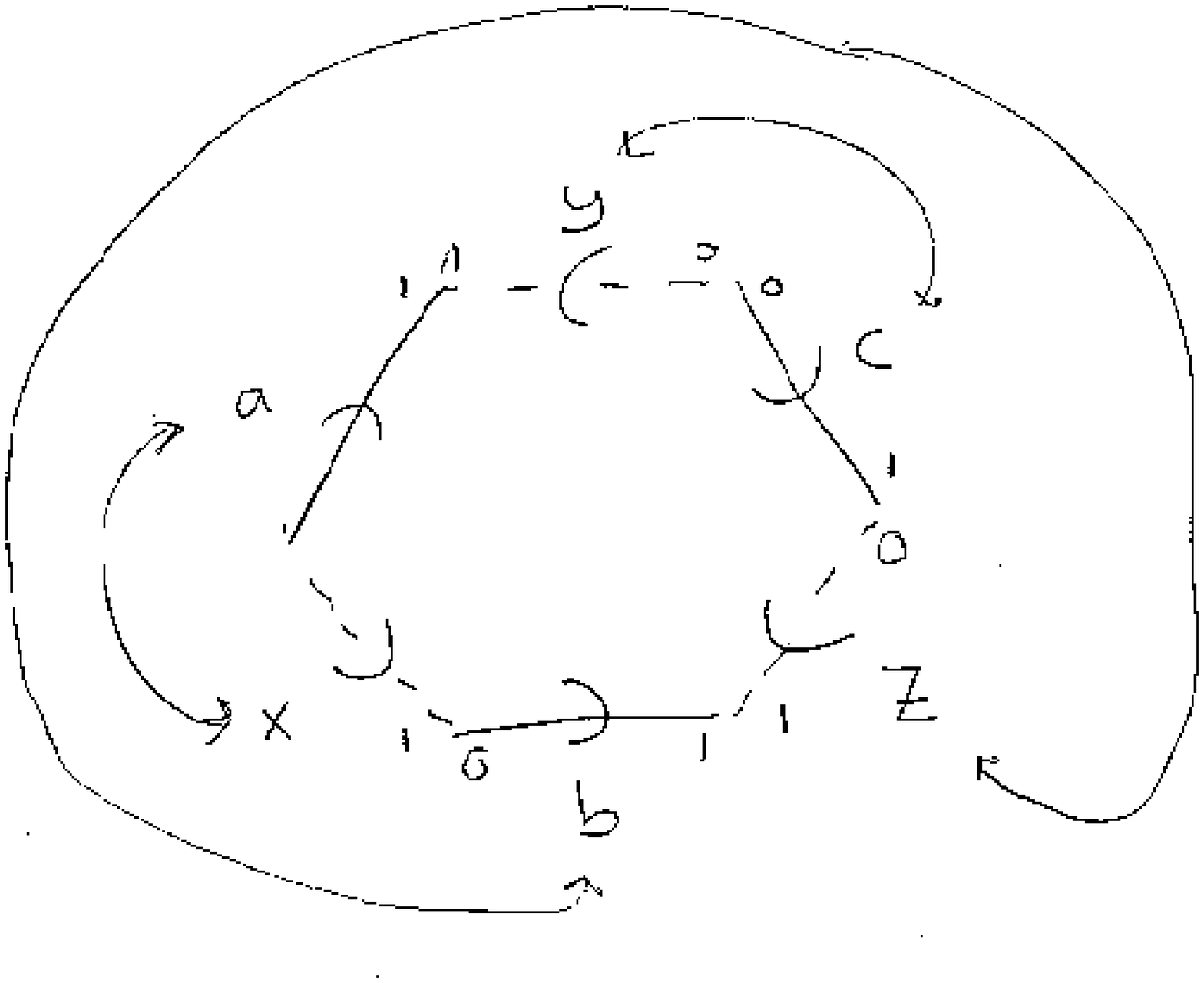}{fig:parens}
{
The parentheses match as indicated.
The curve indicating the matching between
$b$ and $z$ goes the long way around because that's how you go to find
the matching between them.
}
This is a minor change, but it is very suggestive,
and what it suggests is that we match parentheses.
In this case $a$ matches $x$, $y$ matches $c$, and $b$ matches $z$.
Thus in this case we can take as our correspondence
$g(a) = x$, $g(b) = x$, $g(c) = y$.

It is important to realize here
that the notion of parenthesis-matching is perfectly definite,
in the sense that there can never be more than one parenthesis to match
any given one,
and the relation of matching is symmetric (I match you if and only if you
match me),
and matching parentheses always have opposite colors.
All of these properties follow from the observation that
to find the parenthesis that matches a given starting parenthesis,
you can use the familiar counting trick,
where you count 1 for the starting parenthesis and then go off along
the string of parentheses in the direction of the inside of the starting
parenthesis
(that is, {\em backwards} along the corresponding arrow in our original model),
increasing the count whenever you encounter an open parenthesis,
and decreasing the count whenever you encounter a close parenthesis.
The parenthesis matching the starting parenthesis is the first parenthesis
for which the count is 0.
Thus in our example to find the parenthesis that matches $c$
you count 1 for $c$ and 0 for $y$, so $y$ matches $c$.
To find the parenthesis matching $b$ we count 1 for $b$,
2 for $x$,
1 for $a$,
2 for $y$,
1 for $c$,
0 for $z$,
so $z$ matches $b$.

Of course this matching parenthesis idea doesn't always work,
as is illustrated in Figure \ref{fig:unmatched}.
\fig{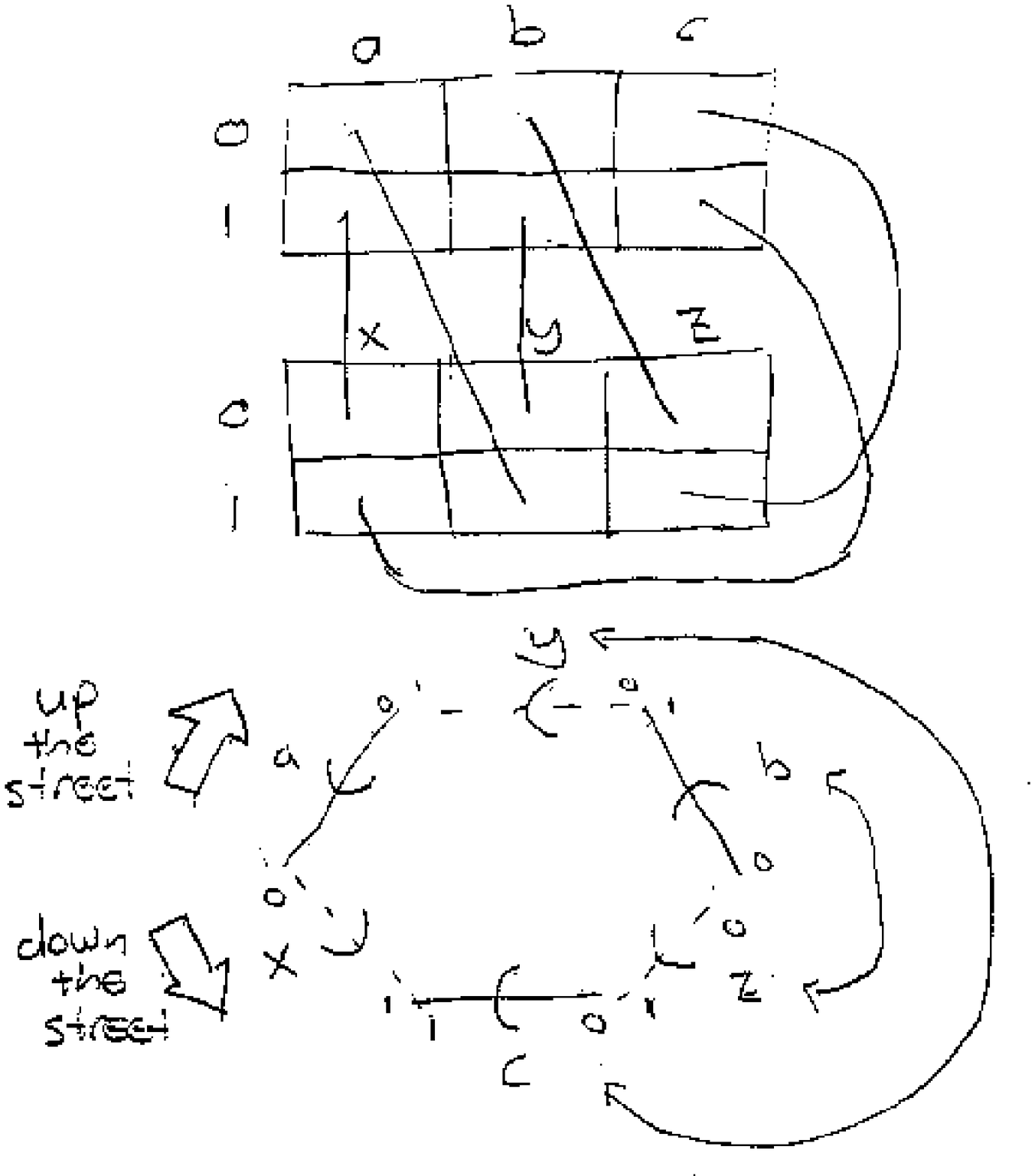}{fig:unmatched}
{
In this example $a$ and $x$ remain unmatched.
They point in the same direction,
and thus determine which way to orient the street.
}
Some of the parentheses may remain unmatched, as
is the case for $a$ and $x$ in this example.
But on a necklace, all the unmatched parentheses point the same way around
the loop,
and we can use this to determine a direction around the loop,
and so in this case we can go back to our idea of having every blue parenthesis
marry the red parenthesis just up the street!

So if there are only necklaces, and no infinite strings, we proceed as follows.
We treat each necklace separately.
First we try matching parentheses.
If this, works, fine.
If not, all the unmatched parentheses face the same way around the loop,
and the way they face determines a direction around the loop.
Let us say that the direction in which the arrows point is `down the street';
that is, the open ends of the unmatched parentheses face `up the street'.
The unmatched parentheses alternate in color,
so if every blue parenthesis marries the next
unmatched (red) parenthesis
up the street,
we will have successfully paired up all the parentheses on the street.

There is another approach that we could have taken to deal with the case
when some of the parentheses are unmatched,
and that is to use the direction in which the unmatched parentheses point
to determine which direction around the loop is `up the street',
and once this has been determined,
to forget about the parenthesis-matching idea altogether, and have
every blue parenthesis, whether matched or not, marry the
red parenthesis next-door-up-the-street.
The result of this second scheme for the example above is illustrated
in Figure \ref{fig:rematch}.
\fig{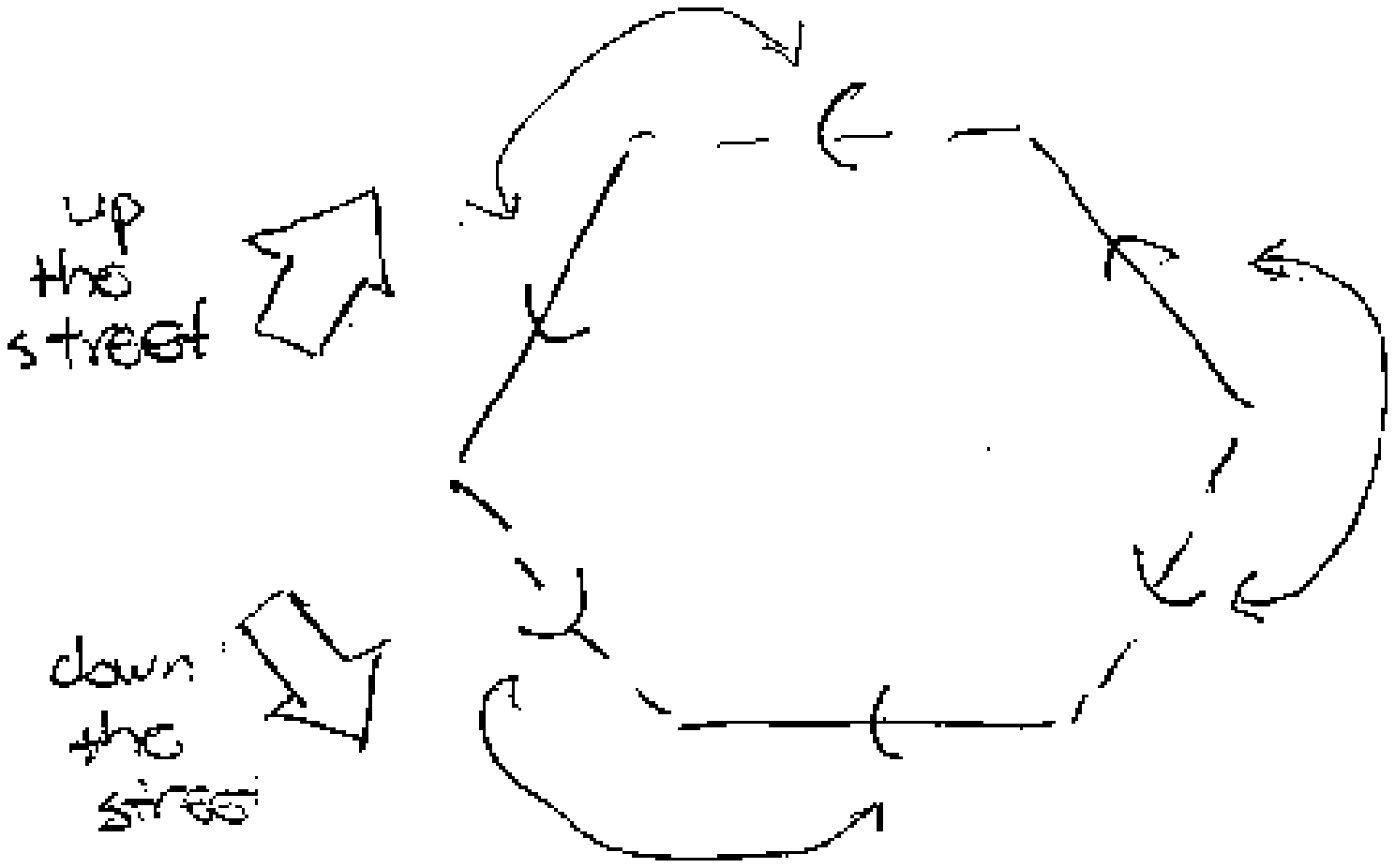}{fig:rematch}
{
According to this strategy,
once it has been determined which direction is `up the street',
the parenthesis-matching idea is abandoned and the blue parentheses
marry next-door-up-the-street,
while the red parentheses marry next-door-down-the-street.
}
This second idea seems a little cruel:
These pairs of matching parentheses think they're meant for each other,
but then we call the whole thing off and tell them to marry someone
else.
But in the infinite case we may have to do this anyway, as we are about
to see.

When $A$ and $B$ are infinite, then there are likely to be infinite strings
of arrows in addition to finite necklaces.
We can still try matching parentheses on these infinite strings, but
when we're done,
the unmatched parentheses that are left over may not all face the same
direction:
There can be two lots that stand back to back.
(See Figure \ref{fig:lots}.)
\fig{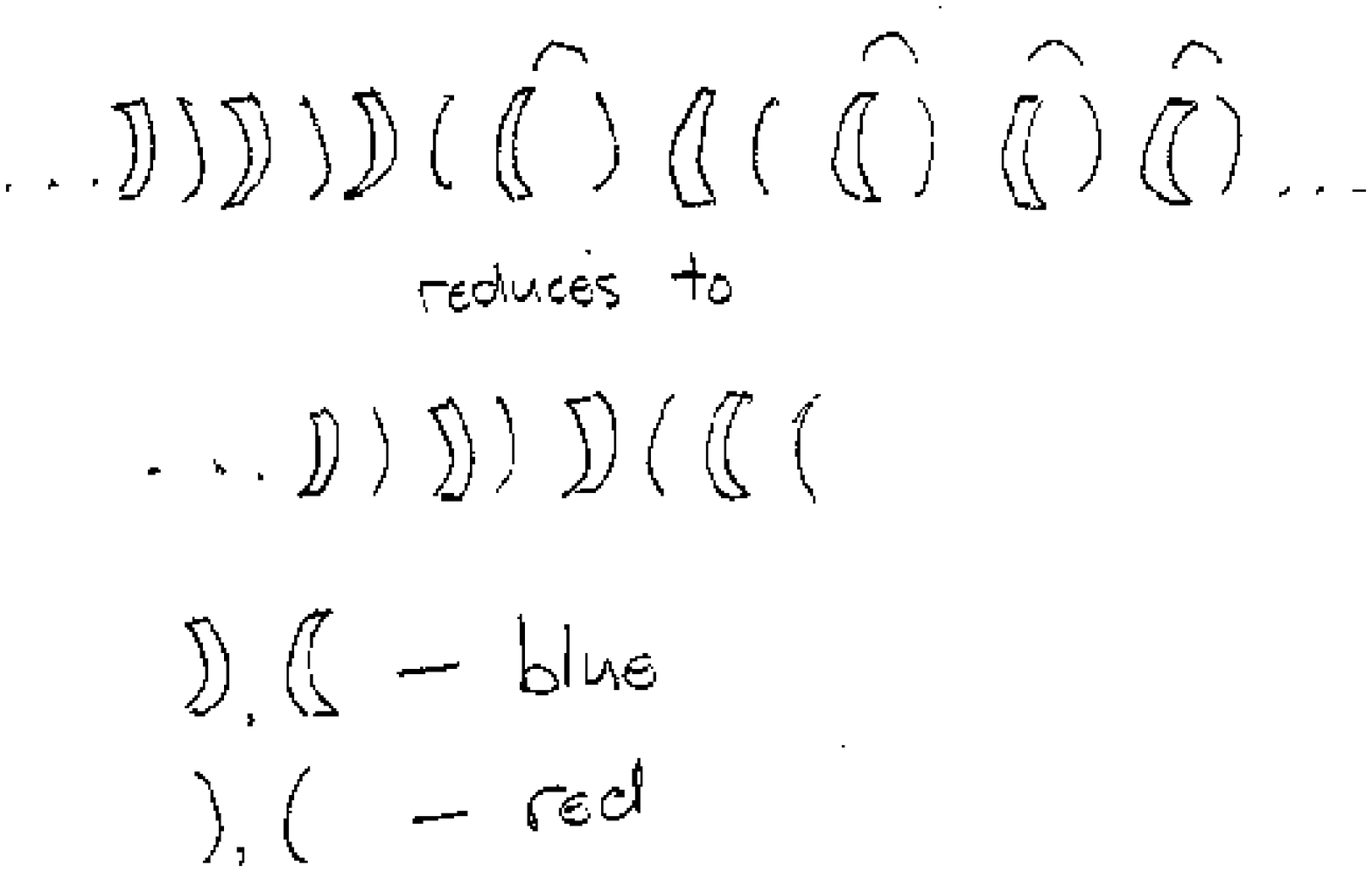}{fig:lots}
{
After removing matching parentheses from an infinite chain,
we can be left with two lots of parentheses facing away from each other.
In this example, the leftward lot is infinite,
while the rightward lot is finite.
}
And while the unmatched parentheses will still alternate in color,
there may be more of one color left than there are of the other;
for example, there might well be only one parenthesis that goes unmatched.
Some of the possibilities are illustrated in Figure \ref{fig:possibilities}.
\fig{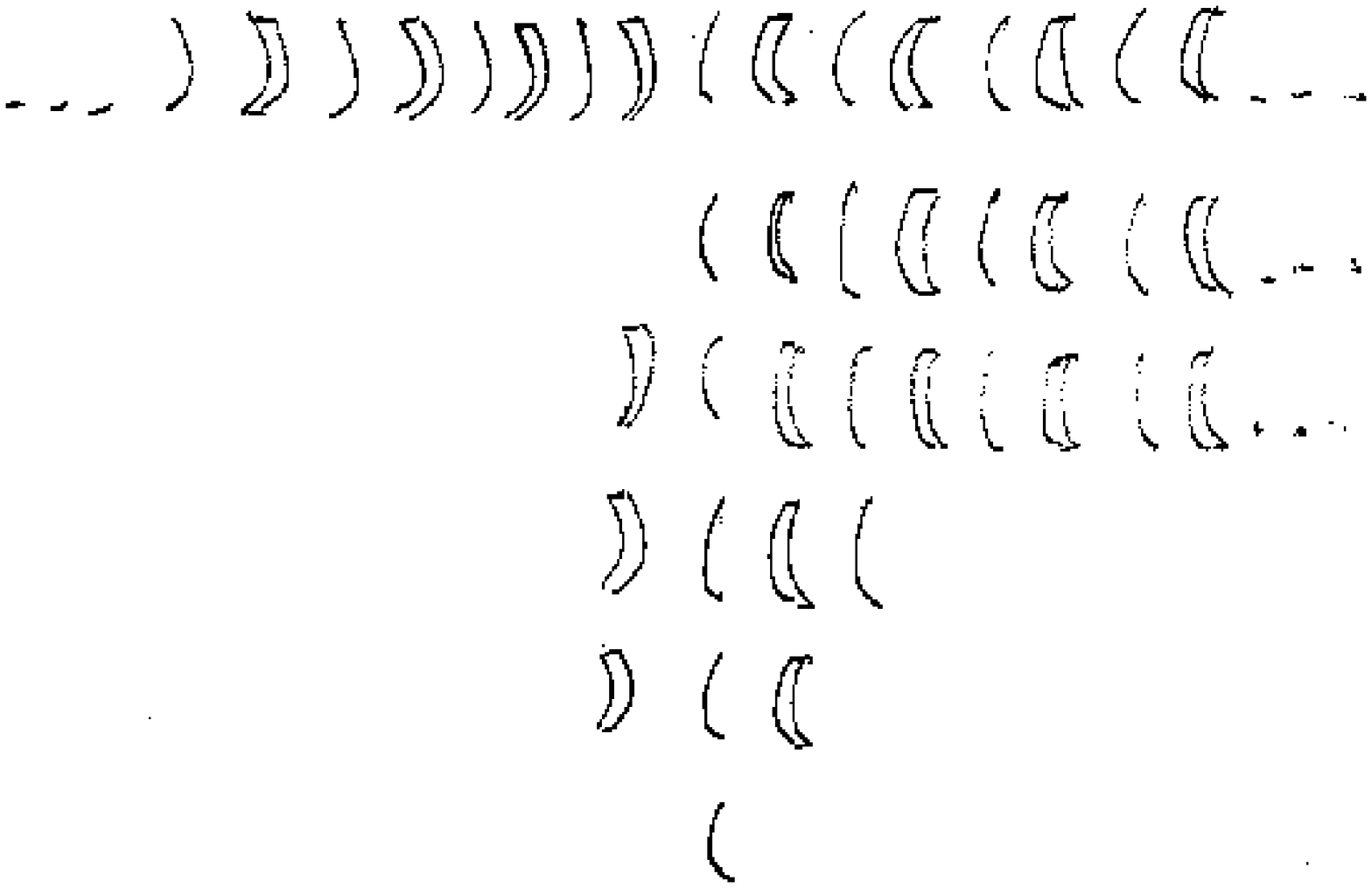}{fig:possibilities}
{
Various possibilities for left-over parentheses:
two infinite lots;
a single infinite lot;
an infinite lot and a finite lot;
an even number;
an odd number;
a single parenthesis.
}

Now if there are no unmatched parentheses on the string we're considering,
then we're in fine shape, because we've paired everybody up.
And if there are some
unmatched parentheses but they all face the same way, then we're still
in fine shape, because this allows us to single out a direction `up the
street', which we can use to implement the plan of having all the
blue parentheses marry the red parenthesis next-door-up-the-street.
Note that in this case, we are not even going to consider letting the
matching parentheses remain together,
because there may be more blue
than red parentheses that are unmatched.

This leaves the case where not all the unmatched parentheses point the
same way.
But as we observed above, in this case the parentheses come in two lots.
Within each lot the parentheses face the same way, namely away from the
other lot.
Thus of the unmatched parentheses, there are two that stand back to back
(that is, are far as unmatched parentheses go; there may be matched
parentheses that intervene).
And since the colors of the unmatched parentheses alternate,
one of the back-to-back pair is blue and the other is red,
so the description `the blue unmatched parenthesis that stands back to
back with a red unmatched parenthesis'
singles out a particular parenthesis of the string.
Of course any time we can single out one parenthesis,
we can determine a direction along the string,
namely, the direction that the arrow corresponding to that parenthesis
points (the direction opposite to the direction this parenthesis faces).
So once again we are in good shape.

To recapitulate:  On each string, we begin by matching parentheses.
If all parentheses are matched, we're all set.  If not, then by examining
the unmatched parentheses we can determine a direction along the string,
and use it to implement the marry-the-parenthesis-next-door strategy.
So we can divide by two.

\section{Division by three}

We are given a bijection
$f:\{0,1,2\} \cross A \rightarrow \{0,1,2\} \cross B$,
and we want to produce a bijection
$g:A \rightarrow B$.
As usual, we assume that $A$ and $B$ are disjoint.

We can visualize the situation by imagining that the set $A$ consists of
a bunch of blue plastic
$30\degree$-$60\degree$-$90\degree$ triangles,
and the set $B$ of red plastic
$30\degree$-$60\degree$-$90\degree$ triangles.
We label the vertices of each triangle 0,1,2 in order of increasing angle,
and we interpret the bijection $f$ as a set of instructions for joining
the triangles together at the vertices.
(We may be called upon to bend and stretch the triangles in order to fit
them together, but of course this is irrelevant to what we're concerned
with here.)
As in dividing by two,
after attaching the triangles at the vertices,
the triangles form independent connected components,
some of which may be finite and some infinite.
An example is illustrated in Figure \ref{fig:triangles}.
(This example looks like an infinite tree homogeneous of degree three,
and in fact it is sufficient to consider this case
`because everything is canonical so you can work up in the universal
cover',
but we won't need to make use of this fact, so
never mind.)
\fig{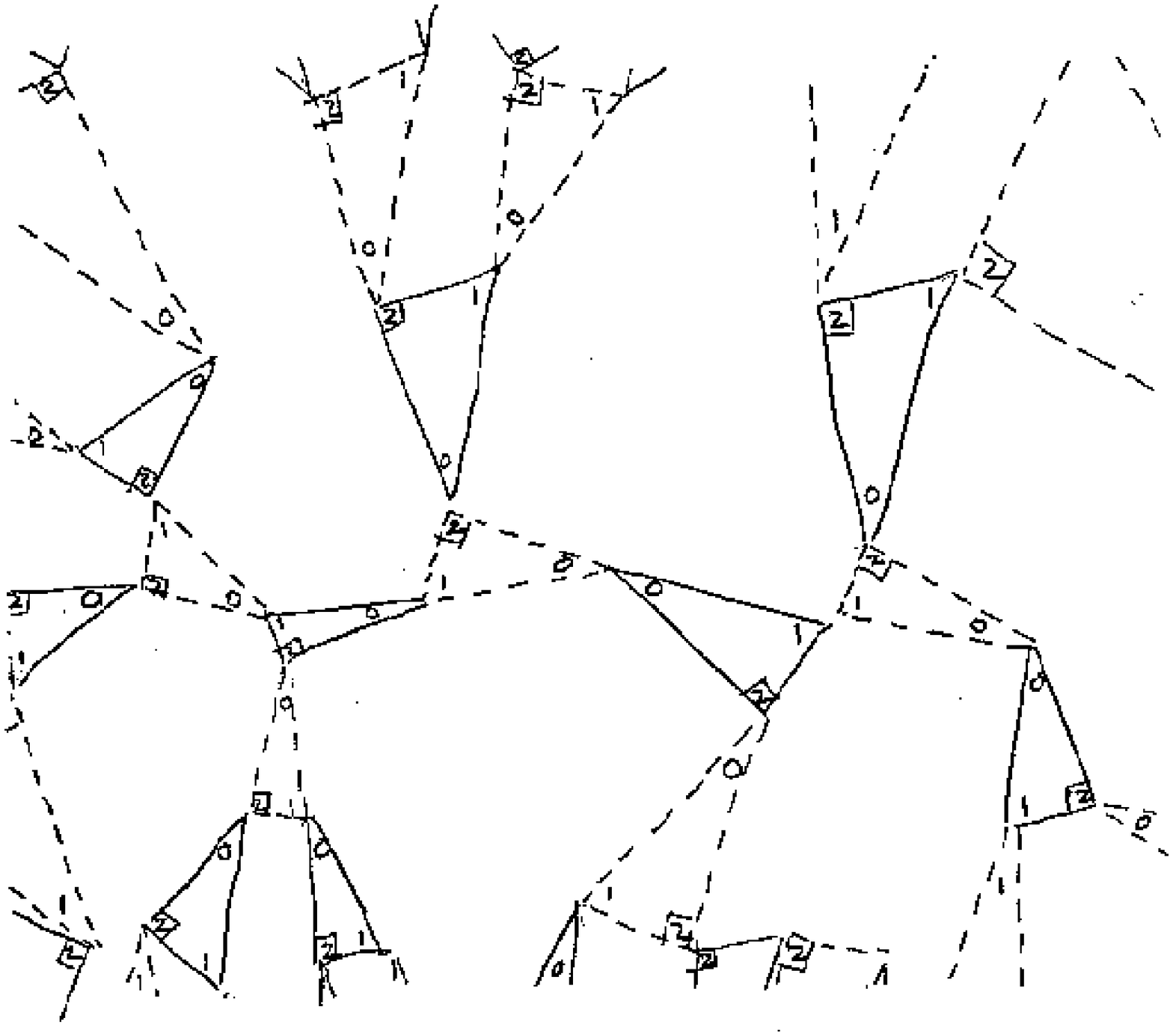}{fig:triangles}
{
Red and blue triangles belonging to a single component.
This example is special in that it is simply connected:
There are no interesting closed loops along the edges of the triangles.
}

Now we tried a lot of ideas in our attempts to adapt the argument that
worked in dividing by two to the case at hand.
But time after time,
we found that despite all our precautions,
every scheme we could come up with for pairing up the triangles
belonging to a given component
had the
defect that we could arrive at a situation where there were an infinite
number of unpaired triangles left over that were all of the same color.
We could see how to handle the case where  there were only a finite
number of monochromatic unpaired triangles,
but we just couldn't handle an infinite number.
So finally we turned in despair to Tarski's article
\cite{tarski:cancellation}
and found,
as noted in the Introduction above,
that as far as Tarski could recall,
Lindenbaum's proof had involved,
in addition to a Sierpi\'{n}skian construction
(the kind of thing we'd been trying so hard to make work),
an appeal to the following lemma, due to Tarski.

{\bf Tarski's Lemma.}
If $A \setgeq B$ and $3 \cross A \setleq 3 \cross B$ then $A \seteq B$.

But of course this lemma is precisely what we needed to finish our argument!
Because say that after doing our utmost to pair up the triangles on a
certain component,
we find to our dismay that we've paired up all the red triangles with
blue triangles,
and we have some unpaired blue triangles left over.
Think of $A$ here as corresponding not to the
original set $A$ of blue triangles,
but to the set of blue triangles belonging to the component of triangles
we're looking at,
and similarly for $B$.
Our incomplete pairing entails the inequality $A \setgeq B$,
because for each red triangle we've found a unique blue triangle,
and the inequality $3 \cross A \setleq 3 \cross B$
follows from the equality $3 \cross A \seteq 3 \cross B$,
which expressed the fact that we have a one-to-one correspondence
between the vertices of the blue triangles in the component and those of
the red triangles.
So Tarski's lemma is going to guarantee that we can convert our
partial pairing into a one-to-one correspondence between the red and blue
triangles belonging to the component we're considering.

Now we don't propose to say anything much about the kind of fiddling needed
to produce on a component a pairing that uses up either all the red or blue
triangles.
This takes a certain amount of thought and work,
but we are confident that you, the reader,
can find such a procedure for yourself,
if you are inclined to do so.
What we do want to explain is a proof of Tarski's lemma,
since this is the missing ingredient that you'll need to convert
a partial pairing on a component into a one-to-one correspondence.
Before we can explain the proof,
however,
we need to explain how to subtract.

\section{Subtraction}

You've probably been wondering why we haven't discussed subtraction
already:
Why all the fuss about division when it's not even clear how to subtract?

In addition to the notation that has been introduced so far,
we want to add the following:
$\omega = \{0,1,2,\ldots\}$ will denote the set of natural numbers.
$A \plus B$ will denote the set $0 \cross A \union 1 \cross B$.
This notation is meant solely to avoid difficulties that can arise
when you want to look at the union $A \union B$, but you're worried
that $A$ and $B$ might not be disjoint.
We will probably abuse logic and think of $A$ and $B$ as being subsets
of $A \plus B$.

{\bf Subtraction Theorem.}
For any finite set $C$,
$A \plus C \seteq B \plus C \implies A \seteq B$.

{\bf Proof.}
What we must show is that given disjoint sets $A$, $B$, and $C$,
where $C$ is finite, and a bijection $f:A \union C \to B \union C$,
we can determine a bijection $F:A \to B$.
To do so, we can visualize the situation much as we did in the proof
of the Cantor-Schr\"{o}der-Bernstein theorem.
We make a copy $C'$ of $C$ disjoint from $A$, $B$, and $C$,
and convert $f$ into a bijection
$f':A \union C \to B \union C'$.
Now we form the same kind of two-colored directed graph as before,
where the vertices of $A \union C$ are blue,
the vertices of $B \union C'$ are red,
blue arcs run from $x \in A \union C$
to $f'(x) \in B \union C'$,
and for each $y \in C$ a red arc runs from the copy of $y$ in $C'$
to $y$.
(See Figure \ref{fig:cobordism}.)
The bijection $F$ is now staring us in the face:
to get from $x \in A$ to $F(x) \in B$,
we just follow the arrows.
Thus subtraction is possible. $\qed$
\fig{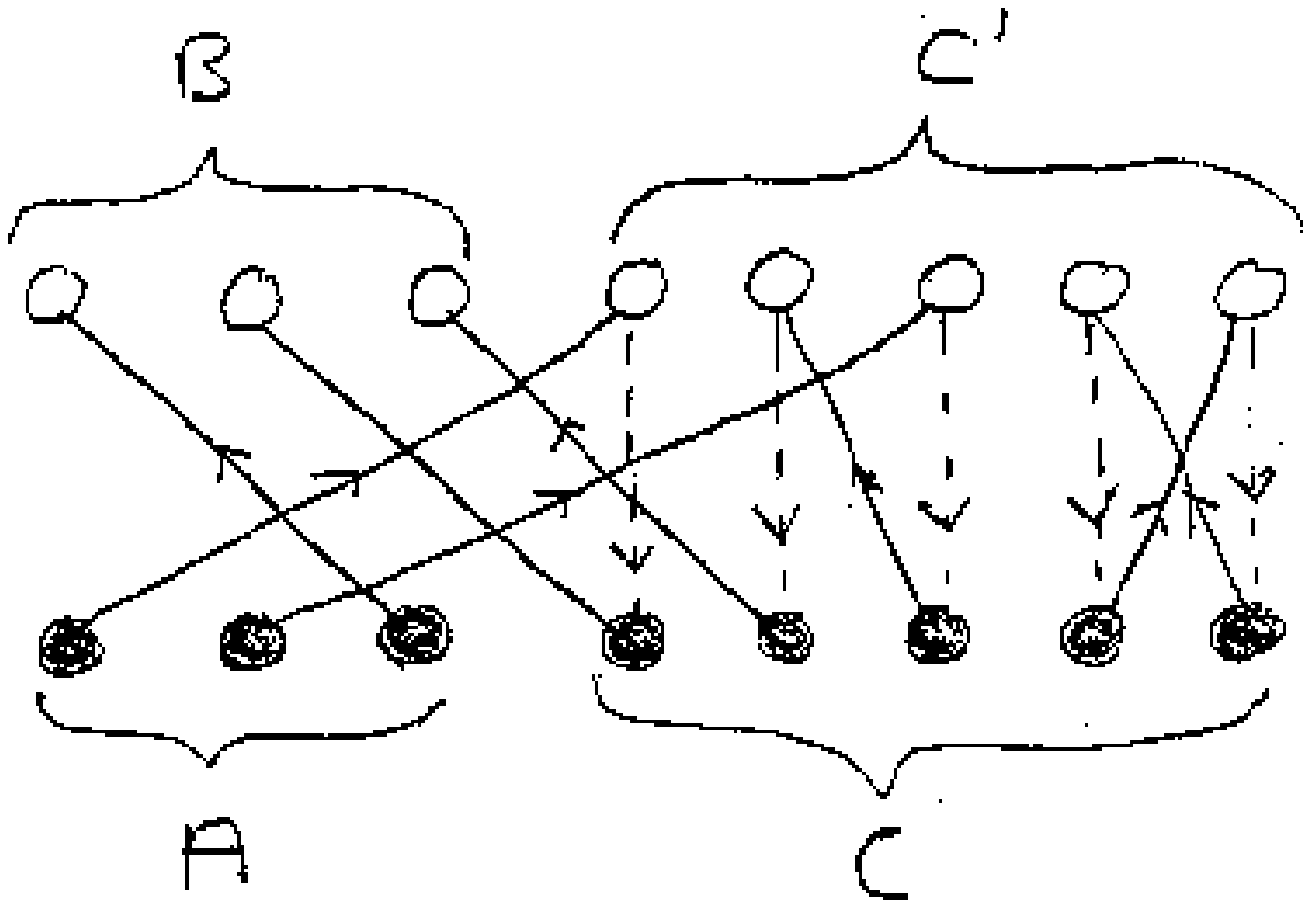}{fig:cobordism}
{
Subtraction,
illustrated in this figure,
is as simple as following the arrows.
It's really just 0-dimensional cobordism theory in disguise.
}

{\bf Remark 1.}
Why does this method of {\em subtracting} involve {\em adding} arrows to
the diagram of $f'$?
Because the arrows run {\em backwards} from $C'$ to $C$.

{\bf Remark 2.}
In dynamical terms,
we can describe the bijection $F(x)$ as obtained by iterating $f$ starting
with $x$
until the image lies in $B$.
This makes sense because if $f(x)$ does not lie in $B$ then it lies
in $C$, so we can apply $f$ to it, and if $f(f(x))$ does not lie in
$B$ then it lies in $C$, so we can apply $f$ to it, etc.
Furthermore, since $f$ is an injection, we can never encounter the
same element of $C$ more than once in this sequence of iterates of $x$,
and since $C$ is finite we must finally emerge in $B$.

{\bf Remark 3.}
In the foregoing proof we have relied on the fact that $C$ is finite,
and indeed things can go wrong if $C$ is infinite,
as illustrated in Figure \ref{fig:lost}.
\fig{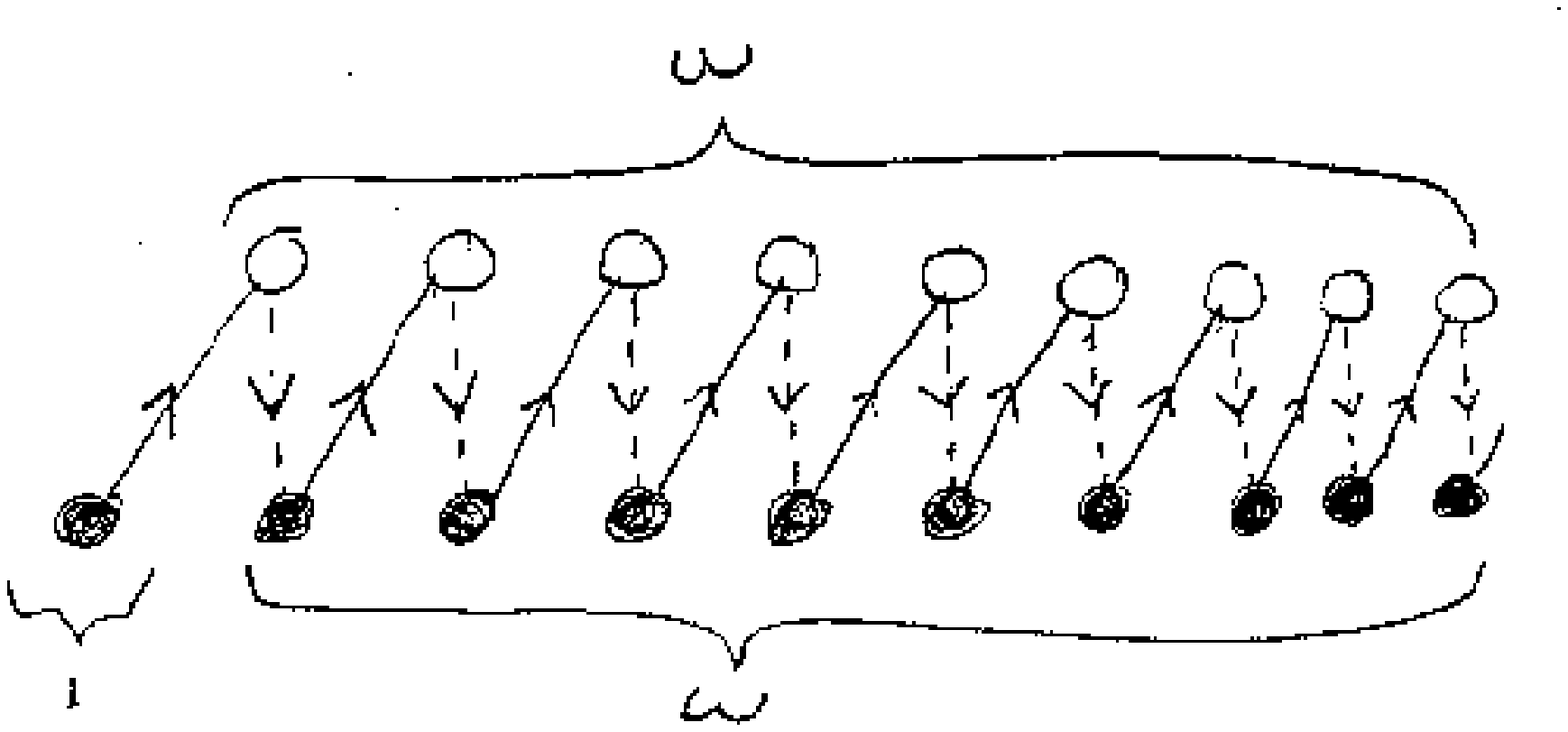}{fig:lost}
{
This is what happens when you try to use subtraction to show that
$1 \plus \omega \seteq 0 \plus \omega \implies 1 \seteq 0$.
}
This figure illustrates the fact that even though
$1 \plus \omega \seteq 0 \plus \omega$, it is not true that $1 \seteq 0$.
In this case,
when we follow the path starting at the unique element of the one-element set
$A$,
we find that it bounces back and forth between the two copies of $C$ and never
emerges.
A set 
is said to be {\em Dedekind-infinite} if there is an injection mapping
the set properly inside itself.
A closer look at the proof above reveals that a set $C$ can be
cancelled (subtracted off)
if and only if it is Dedekind-finite.
This is stronger than the result stated because
in the absense of the axiom of choice,
there may be Dedekind-finite sets that are `infinite'
in the sense that they can't be mapped bijectively to $\{0,1,\ldots,n-1\}$
for any $n$.

\section{Swallowing}
If $A \plus B \seteq A$
we will say that $A$ {\em swallows} $B$,
and write $A \setgg B$,
bearing in mind the danger inherent in the fact that very often
$A \setgg A$,
e.g. $ \omega \setgg \omega$.
Note that a set is Dedekind-infinite
if and only if it swallows some non-empty set.

While we have defined $A \setgg B$ to mean
$A \plus B \seteq A$,
we could equally well have made the condition be that $A \plus B \setleq A$,
because the obvious injection shows that $A \setleq A \plus B$,
and thus by Cantor-Schr\"{o}der-Bernstein
$A \plus B \setleq A \implies A \plus B \seteq A$.
In some cases, we will find it convenient to produce an injection from
$A \plus B$ to $A$, and then let
Cantor-Schr\"{o}der-Bernstein do the work of turning this injection into
a bijection.

The phenomenon of swallowing can be visualized by means of the
{\em Hilbert Hotel} model.
(See Figure \ref{fig:hotel}.)
\fig{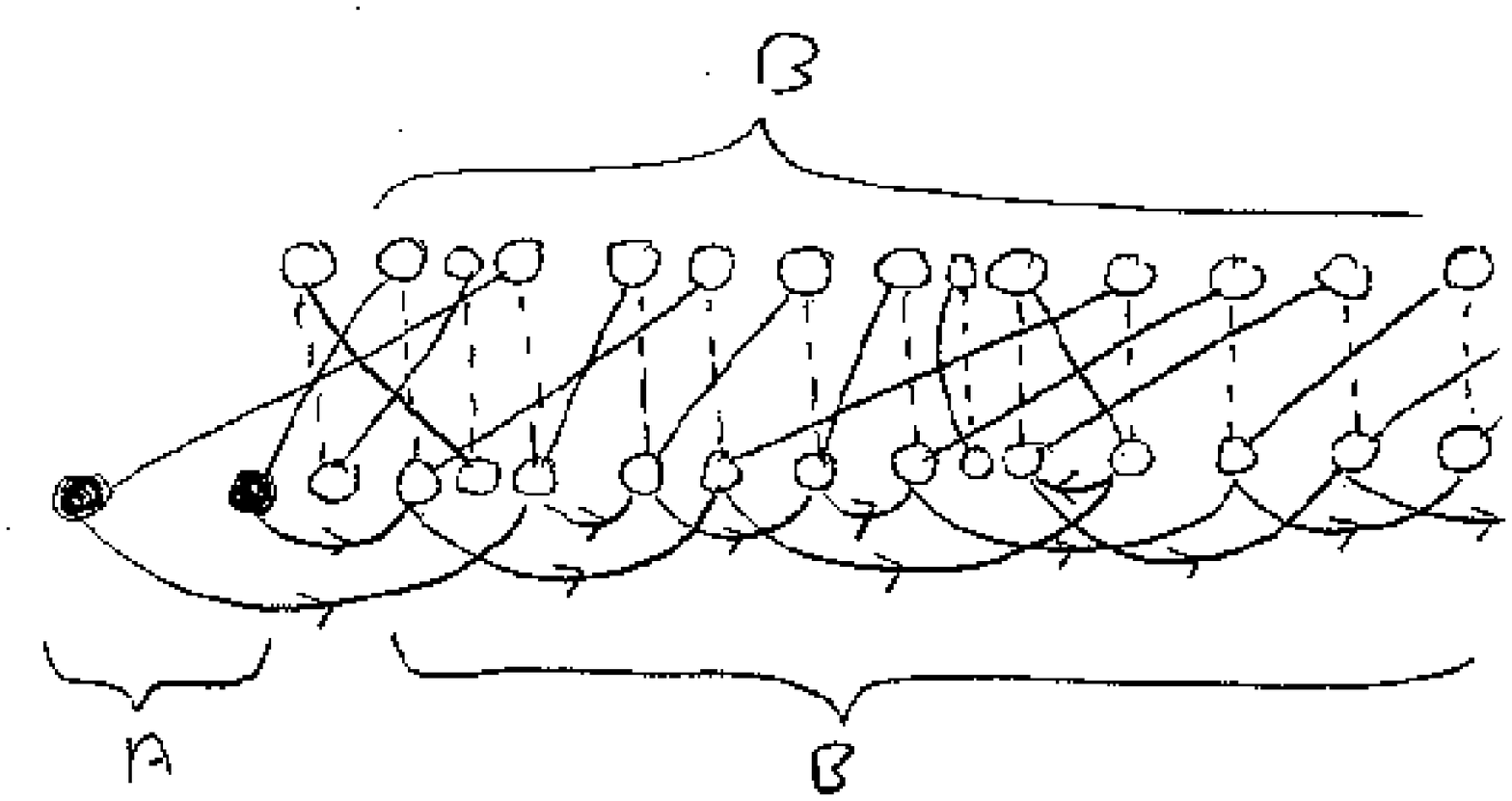}{fig:hotel}
{
A bijection between $A \union B$ and $B$ yields a collection of disjoint
sequences in $B$,
one for each element of $A$.
}
We think of the elements  of $A$ as representing
the rooms in a hotel that we are managing.
The hotel is chock full, but then a whole bunch of new guests
come,
represented by the elements of $B$.
``No problem,'' we tell them,
``we'll fit you all in.''
To make room for $x \in B$,
we tell the guest in $f(x)$ to move to room $f(f(x))$;
we tell the guest in room $f(f(x))$ to move to room $f(f(f(x)))$,
etc.
This model is also called the {\em Roach Hotel} model,
after a famous cockroach trap called the Roach Hotel,
where ``Roaches check in---but they don't check out!''

{\bf Theorem.}
$B \setll A$ if and only if $B \cross \omega \setleq A$.

{\bf Proof.}
For each element of $x \in B$ we get
a whole $\omega$'s-worth of elements of $A$ by
taking the orbit of $x$ under the injection $f:A \plus B \to A$,
namely the sequence of rooms through which
guests are moved to make way for $x$.
$\qed$

\section{Proof of Tarski's lemma}

Now we come to the proof of Tarski's lemma
that if $A \setgeq B$ and $3 \cross A \setleq 3 \cross B$
then $A \seteq B$.
Because the proof is broken down into a series of lemmas,
there is a danger that we will lose sight of what is actually going on,
and what was supposed to be a proof will turn out to be nothing
more than a verification.
To try to avoid this pitfall,
when we're finished with all the lemmas
we will go back and try to really understand Tarski's lemma
by visualizing how it goes when you put it all together.

{\bf Lemma 1.}
If $A \setll C$ and $B \setll C$ then $A \plus B \setll C$.

{\bf Proof.}
First one bunch of guests check in, and then another bunch of guests check in.
$\qed$

{\bf Lemma 2.}
If $A \setll B_1  \plus  B_2  \plus  B_3$ then we can write
$A = A_1  \plus  A_2  \plus  A_3$
where $A_1 \setll B_1$, $A_2 \setll B_2$, $A_3 \setll B_3$.

{\bf Proof.} 
We imagine a bunch of guests checking into a hotel consisting of three
separate buildings.
Each guest causes a cascade of guests along a sequence of rooms of the hotel.
(See Figure \ref{fig:threebuildings}.)
\fig{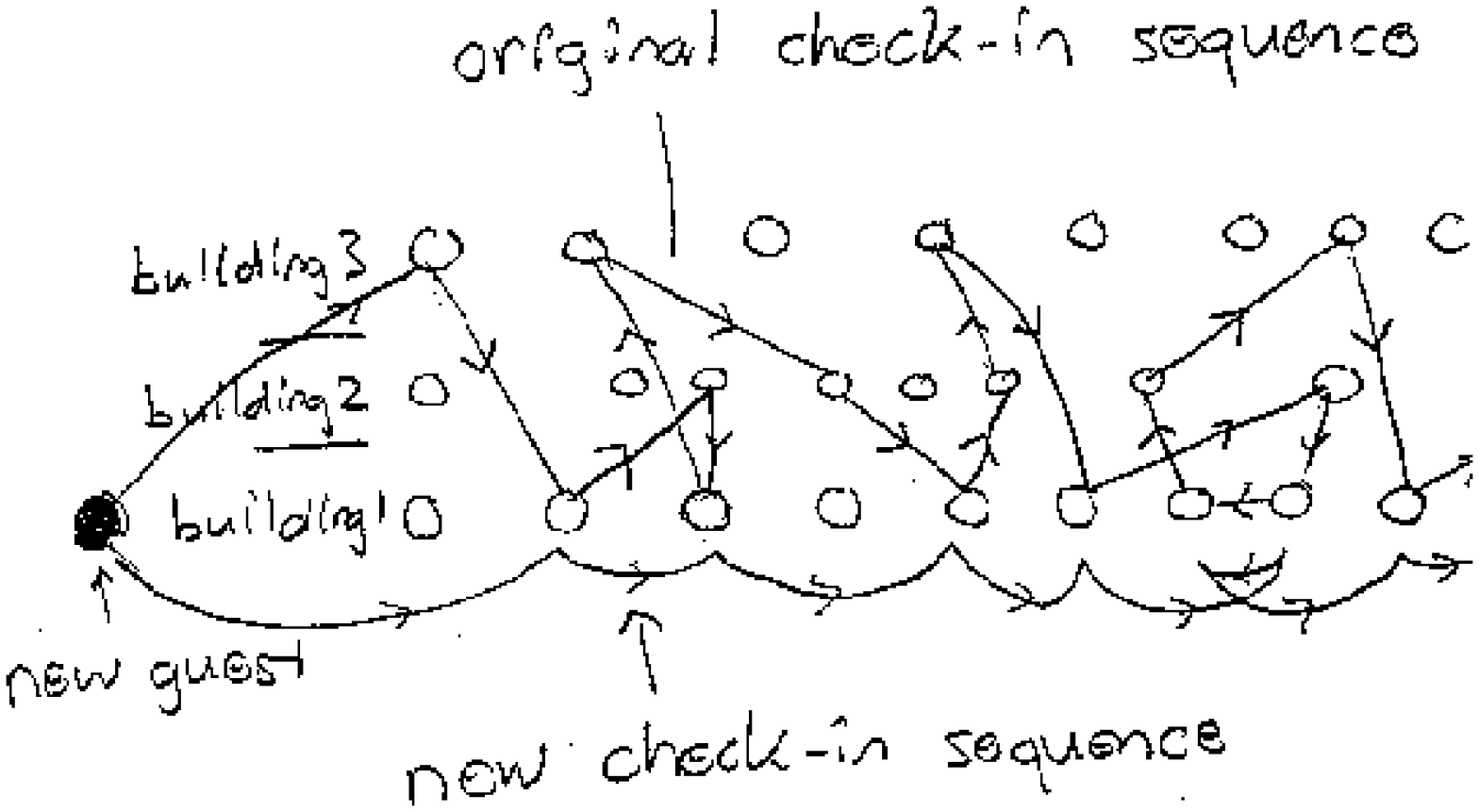}{fig:threebuildings}
{
It appears that we can check this new guest in without disturbing
any of the guests who are currently in buildings 2 and 3.
}
For every guest, the corresponding sequence involves an infinite number of
rooms in at least one of the three buildings.
If it involves an infinite number of rooms
in building 1, we put this new guest into the first room in building 1,
move the guest in that room into the next room in the sequence that belongs
to building 1, etc.
Otherwise if the sequence involves an infinite number of rooms from building 2,
we check the new guest into the first room in the sequence in building 2, etc.
Otherwise we use the rooms in building 3.
and similarly for the sequences involved in checking guests into buildings
2 and 3.
Note that it is important that
for a given building,
the sequences of rooms involved in checking
guests into that building
under this modified scheme are all disjoint.
$\qed$

{\bf Lemma 3.}
If $A \setll 3 \cross B$ then $A \setll B$.

{\bf Proof.}
Combine Lemmas 1 and 2.
$\qed$

{\bf Tarski's Lemma.}
If $A \setgeq B$ and $ 3 \cross A \setleq 3 \cross B$ then $A \seteq B$.

{\bf Proof.}
If $A \setgeq B$ then $A \seteq B \plus C$ for some $C$
(namely, the complement of the image of $B$ under some injection from $B$
into $A$).
We want $A \seteq B$, i.e.\ we want $C \setll B$.
We have
\begin{eqnarray*}
&&3 \cross (B \plus C) \setleq 3 \cross b 
\\&\implies&
3 \cross B  \plus  3 \cross C \setleq 3 \cross B
\\&\implies&
3 \cross C \setll 3 \cross B
\\&\implies&
C \setll 3 \cross B
\\&\implies&
C \setll B
\dispqed
\end{eqnarray*}

Well there you have it.
Now let's take a closer look at it,
and see if we can see it all happening at once.
So suppose that $3 \cross A \setleq 3 \cross B$,
or rather,
to simplify matters, suppose that $3 \cross A \seteq 3 \cross B$.
We construct our graph out of blue and red triangles as before,
only now instead of glueing pairs of corresponding vertices together,
let's connect them by bits of white string.
We're given a matching between the red triangles and some, but not all,
of the blue triangles,
and we want to convert this matching into a one-to-one correspondence
by finding mates for the blue triangles that have been left out in the cold.
Let's connect each red vertex to the similarly-labelled vertex of the
corresponding blue triangle by a bit of black string.
Note that a blue vertex has a black string attached to it
if and only if it has been matched to some red triangle.
Note also that the black strings may have joined together triangles
that belonged to separate white-string components,
because we're not assuming that the partial pairing between red and blue
triangles respects the components of the white-string graph,
though this will be the case in the application of Tarski's lemma
to division by three.

Starting from vertex 0 (say) of some unmatched blue triangle $x$,
there is a unique infinite path along the strings with colors
alternating white, black, white, black, \ldots.
As we follow this path,
the vertices we encounter belong to triangles colored
blue, red, blue, red, \ldots.
Every time we finish traversing one of the white strings we find ourselves
at a vertex of some red triangle;
every time we finish traversing one of the black strings we find ourselves
at a vertex of some blue triangle.
If infinitely many of these blue vertices are labelled 0,
we call the triangle $x$ a {\em 0-triangle},
and associate to it the infinite sequence of blue triangles whose
0-vertices we encountered on our walk.
All the sequences belonging to 0-triangles are disjoint,
and each triangle in the sequence other than the starting triangle $x$
has a partner in the original matching.
This means that we can rearrange the matching by having each blue triangle
pass its partner one notch up the sequence,
so that now all the blue triangles of the sequence including $x$ will
have partners.
{\em Before doing this},
however,
we divide the remaining unpaired blue triangles into 1-triangles
and 2-triangles,
and associate to each the appropriate sequence of blue triangles.
After partners have been passed up the 0-triangle sequences,
the {\em new} partners are passed up the 1-triangle sequences,
and then the {\em new, new} partners are passed up the 2-triangle sequences.
Get it?

In this description,
we assumed for simplicity that we had $3 \cross A \seteq 3 \cross B$
in place of $3 \cross A \setleq 3 \cross B$.
Why?
We could still construct the complex of triangles,
with the only difference being that now some red triangles might have
free vertices.
This wouldn't have affected the argument we've just given.
We will go into this more closely in the next section,
where we consider what has to be done
to extend division by three to inequalities.

\section{Dividing an inequality by three}

It behooves us to consider whether the proof that we have
given of the fact that $3 \cross A \seteq 3 \cross B$ implies $A \seteq B$
can be adapted to
show that $3 \cross A \setleq 3 \cross B$ implies $ A \setleq B$.
That will depend to some extent on you, the reader,
because we have been relying on you to supply a procedure which,
on each connected component of the complex of triangles,
will produce a partial matching such that the triangles that
remain unmatched all have the same color.
As indicated at the end of the last section,
we now want to consider the complex of triangles associated to an
injection of $3 \cross A$ into $3 \cross B$,
a complex in which some of the red vertices may remain free.
The question is,
can your method be adapted to handle this case?
Actually, we are pretty confident that if it can't,
you will be able to come up with a new method that can be made
to work in this more general case.
But now we're going to reveal one of our methods,
a method that works in this general case as well.
This means that when we said above that we were going to leave it
entirely up to you to supply your own method for constructing
partial matchings,
we were lying.
The reason for this is that we felt it would really be much better
for you to find your own method,
instead of reading about ours,
and that if we made you think we weren't going to tell you the answer,
then there was a chance that you would look for a method yourself,
but that if we let you know we were going to reveal the method,
then human nature being what it is,
you would probably want to read about our method rather then thinking
for yourselves.
Anyway,
in return for revealing the secret,
we ask that you refrain from telling anyone else that the secret
is actually revealed here
in this innocuous-sounding section.

So here's how it goes.
We look at a connected component of the complex of triangles.
To each pair consisting of a blue triangle $x$ and a red triangle $y$,
there are various paths through the complex leading from $x$ to $y$.
We can describe a path leading from $x$ to $y$ by saying something like,
`You leave the starting triangle via vertex 0,
which is connected to vertex 1 of a red triangle,
whose vertex 0 is connected to vertex 0 of a blue triangle,
whose vertex 2 is connected to vertex 0 of the destination triangle'.
(See Figure \ref{fig:route}.)
\fig{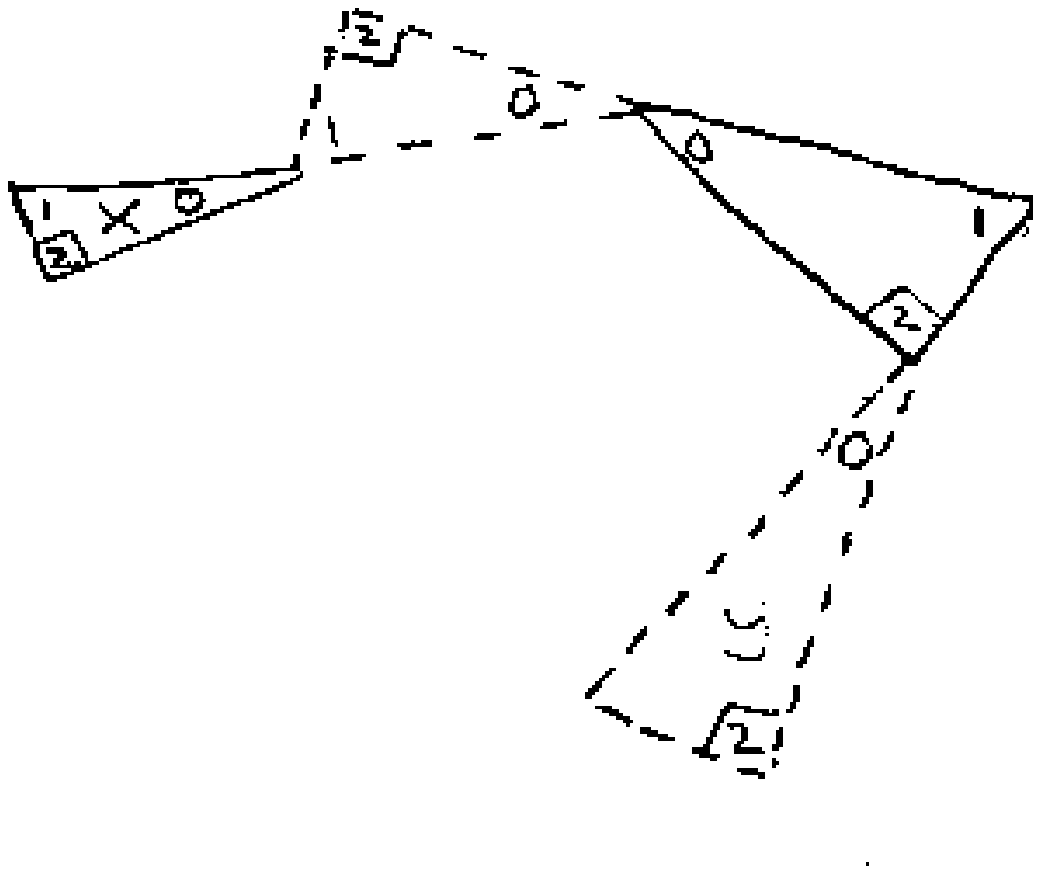}{fig:route}
{
How to get from $x$ to $y$:
`You leave the starting triangle via vertex 0,
which is connected to vertex 1 of a red triangle,
whose vertex 0 is connected to vertex 0 of a blue triangle,
whose vertex 2 is connected to vertex 0 of the destination triangle'.
}
The important aspect of this description is that for a given blue triangle
$x$, there can be at most one red triangle $y$ that you can get to from
$x$ by following these instructions.
We now prescribe an ordering on these descriptions by decreeing that
shorter paths come before longer paths,
and within the finite set of descriptions of paths of a fixed length,
the ordering is---well actually it doesn't really matter what it is,
as long as it's definite.
For example, for the shortest paths, which correspond to triangle
that are adjacent in the complex, 
the description is essentially `Blue vertex $i$ attached to red vertex $j$',
and we may take the ordering of the pairs $(i,j)$ to be lexicographic:
$(0,0) < (0,1) < (0,2) < (1,0) < (1,1) < (1,2) < (2,0) < (2,1) < (2,2)$.
Now we pair up all pairs of blue and red triangles
that are connected by paths
whose description
comes first in our ordering.
In our sample ordering, this means pairs of triangles whose vertices
labelled 0 are stuck together.
Then, we select from among those pairs of blue and red triangles
neither of which has yet been matched up all those pairs whose path
comes next in our ordering of paths.
In our example, this would involve all blue triangles whose vertex labelled
0 is stuck to vertex 1 of a red triangle that has not yet been matched,
i.e.\ whose vertex 0 is not attached to a blue vertex 0.
We keep this up, running through our list of descriptions,
at each stage matching up pairs of
as-yet-unmatched triangles that fit the description.
When we're all done,
which will likely take infinitely long unless we start working faster and
faster as we move down the list,
there may be some blue triangles left over,
or there may be some red triangles left over,
but there can't be both a blue and a red triangle left over,
because these would have been matched when we came to consider the
description of a path connecting them.

If there are red triangles left over, that's fine, because we're trying
to show that $A \setleq B$.
If there are blue triangles left over,
then as far as this component goes we're in the situation where
$A \setgeq B$ and $3 \cross A \setleq 3 \cross B$,
and by now we know very well how to modify the correspondence on this
component so that the blue and red triangles are paired exactly.
So in either case we're OK.

\section{More on division by two}
Of course we can modify the method we used to divide by three
to get a method for dividing by any positive integer $n$,
including $n=2$,
and it is instructive to go back now and compare this new method of dividing
by two with the method we used at the beginning of the paper.
We can interpret the matching of parentheses as resulting from an ordering
on paths,
only now it is not an ordering that extends the partial ordering by length
of path,
because we never pair up adjacent parentheses that face the same direction
while there are non-adjacent matching parentheses left to pair up.
So instead we make all paths that identify matching parentheses come
at the beginning of our ordering,
before any path corresponding to parentheses that face the same
direction.
Evidently, this isn't going to be possible if our ordering is that of a
sequence: first, second, third, \ldots.
Instead we need to allow an ordering whose {\em order type} is more exotic.
In this case,
we can get by with what is called the order type $2 \omega$:
first, second, third, \ldots, infinitieth, infinity plus first,
infinity plus second, \ldots.
If we do this, we can be assured that if the parentheses actually match,
then this is the matching we will choose.
By making sure that shorter paths come before longer paths
once we get beyond descriptions of matching parentheses,
i.e.\ once we get to the infinitieth member of our list of patterns,
we can assure that once all the matching is done, there is only
a single parenthesis left over.
This would simplify things a little,
because we could use this single parenthesis to determine a direction
up the street,
and there would be no call to use Tarski's lemma.
This is related to a phenomenon which we allowed to pass unremarked in
the case of division by three:
If you begin with a simply connected complex, and match and match
and end up with only a finite number of unmatched triangles,
then you can select one of these triangles according to some kind of
ordering rule akin to that used to order paths through the complex,
and once you have a distinguished triangle you can describe explicitly
a matching in a way that depends on your having a well-defined triangle
to refer the description to.
(Go 3 steps north \ldots.)
The whole need for Tarski's lemma stems from the need to deal with
the possibility that you had left over an {\em infinite} number of triangles
of the same color.
This is something that is easy to avoid in dividing by two,
which is apparently what makes division by two so much easier than
division by three.

\section{Division as repeated subtraction}

Our description of matching parentheses in terms of an ordering of paths
that puts paths between matching parentheses before all other paths doesn't
really do justice to the concept.
It is really more natural to describe the process in terms of repeated
subtraction.
You begin by pairing up parentheses that stand face to face.
You take this partial correspondence,
induce in the most natural way the associated correspondence between
sides (inside and outside) of the parentheses involved in the partial
correspondence,
and {\em subtract} this from the original correspondence.
You take the difference correspondence,
and treat it as though it were the original correspondence,
get a new partial correspondence, subtract, \ldots.
You do this infinitely many times,
and if you're lucky,
all the parentheses match and you're finished.
More likely, you have parentheses left over,
and you have to fuss.
The fussing is minimal in the case of division by two,
but for division by three,
it's harder.
In dividing by three,
the analog of matching parentheses is to pair triangles that sit
sharp point to sharp point,
subtract,
pair triangles that sit sharp point to sharp point,
subtract, \ldots.
The first thing you notice is that the complex corresponding to the
difference bijection becomes disconnected as soon as you subtract once,
and once you've subtracted an infinite number of times, you can be
left with all kinds of junk.
In particular, you can be left with an infinite number of lonely blue
triangles, completely disconnected from one another.
This doesn't look so promising, so you try something else.
And so on.

The moral?  {\em There is more to division than repeated subtraction.}

There are, however, interesting cases where there is enough extra
structure around that is preserved by the bijections so that
division reduces to repeated subtraction,
and this leads to interesting results in the field of
{\em bijective combinatorics}.
This will be reported on at length elsewhere,
but it is worth saying a few words about it here to draw the connection
between the finite problem and the infinite problem we just got
through discussing.

In bijective conbinatorics,
the big thing is to produce explicit bijections between sets $A$ and $B$
that you know by some round-about route have the same (finite)
number of elements.
In this case, the bijection should not only be explicit,
in the technical sense we've been using here,
but somewhat pleasing or efficient or natural or something.
Occasionally, you find that if you adjoint the same disjoint set
of gadgets $C$ to $A$ and $B$, then you can describe explicitly
a bijection between $A \union C$ and $B \union C$.
Then you get an explicit bijection between $A$ and $B$ by subtraction.
This technique and similar techniques known collectively as
`the involution principle' have been used to great effect by
Garsia and Milne
\cite{garsiaMilne:rogers},
Remmel
\cite{remmel:partition},
and others.

By the same token,
you can sometimes find a bijection between $A \cross C$ and $B \cross C$,
and if there is enough additional structure around you can convert this
to a bijection between $A$ and $B$ by division, which in this case reduces
to repeated subtraction.
The need to have the additional structure,
so that the division reduces to repeated subtraction,
arises because of the requirement that the correspondence
you end up with be somehow nice.
Without the additional structure,
you can't use repeated subtraction alone,
and things are not nice.

The moral?
{\em There is more to division than repeated subtraction.}

\section{What's wrong with the axiom of choice?}

Part of our aversion to using the axiom of choice stems from our view
that it is probably not `true'.
A theorem of Cohen shows that
the axiom of choice is independent of the other axioms of ZF,
which means that neither it nor its negation can be proved from
the other axioms,
providing that these axioms are consistent.
Thus as far as the rest of the standard axioms are concerned,
there is no way to decide whether the axiom of choice is true or false.
This leads us to think that we had better reject the axiom of choice
on account of Murphy's Law that `if anything can go wrong, it will'.
This is really no more than a personal hunch about the world of sets.
We simply don't believe
that there is a function that assigns to each non-empty set of
real numbers one of its elements.
While you can describe a selection function that will work
for finite sets, closed sets, open sets, analytic sets, and so on,
Cohen's result implies that there is no hope of describing a definite
choice function
that will work for `all' non-empty sets of real numbers,
at least as long as you remain within the world of standard
Zermelo-Fraenkel set theory.
And if you can't describe such a function,
or even prove that it exists without using some relative of the axiom of choice,
what makes you so sure there is such a thing?

Not that we believe there really are any such things as infinite sets,
or that the Zermelo-Fraenkel axioms for set theory are
necessarily even {\em consistent}.
Indeed, we're somewhat doubtful whether large natural numbers
(like $80^{5000}$, or even $2^{200}$) exist in any very real sense,
and we're secretly
hoping that Nelson will succeed in his program for proving
that the usual axioms of arithmetic---and hence also of set theory---are
inconsistent.
(See Nelson
\cite{nelson:predicative}.)
All the more reason, then,
for us to stick with methods which,
because of their concrete, combinatorial nature,
are likely to survive the possible collapse
of set theory as we know it today.
\bibliography{three}
\bibliographystyle{plain}
\end{document}